\documentclass[11pt]{amsproc}

\usepackage{amsmath, amsthm, amscd, amsfonts, amssymb, enumerate, verbatim, newlfont, calc, graphicx, color}
\usepackage[bookmarksnumbered, colorlinks, plainpages]{hyperref}
\usepackage{tikz}

\newtheorem{thm}{Theorem}[section]
\newtheorem{question}[thm]{Question}
\newtheorem{lem}[thm]{Lemma}
\newtheorem{pro}[thm]{Proposition}
\newtheorem{cor}[thm]{Corollary}
\theoremstyle{definition}

\newtheorem{rem}[thm]{Remark}
\newtheorem{defn}[thm]{Definition}

\newcommand{\bQ}{\mathbb{Q}}
\newcommand{\bZ}{\mathbb{Z}}

\usepackage{pstricks}
\usepackage{epsfig}
\usepackage{pst-grad}
\usepackage{pst-plot}

\begin{document}
\author[ I. Al-Ayyoub et. al. 
]{Ibrahim ~Al-Ayyoub$^{1}$,  Mehrdad ~Nasernejad$^{2}$, Kazem Khashyarmanesh$^{2,*}$,  Leslie G. Roberts$^{3}$, and Veronica Crispin Qui$\mathrm{\tilde{n}}$onez$^{4}$}
\title[Results on the normality of  square-free monomial ideals] {Results on the normality of  square-free monomial ideals and cover ideals under some graph operations}
\thanks{$^*$Corresponding author}

\thanks{E-mail addresses: iayyoub@just.edu.jo,
 m$\_$nasernejad@yahoo.com, khashyar@ipm.ir,
  robertsl@queensu.ca, and veronica.crispin@math.uu.se}  
\maketitle

\begin{center}
\emph{
$^1$Department of  Mathematics, 
Sultan Qaboos University, \\
 P.O. Box 31, Al-Khoud 123, Oman\\ 
$^1$Department of  Mathematics and Statistics, 
Jordan University of Science \\
 and Technology, P.O.Box 3030, Irbid 22110, Jordan\\
$^{2}$Department of Pure Mathematics,
 Ferdowsi University of Mashhad,\\
P.O.Box 1159-91775, Mashhad, Iran\\
$^3$Department of Mathematics and Statistics, Queen's University,\\
Kingston, Ontario, Canada, K7L 3N6\\
$^4$Department of Mathematics, Uppsala University,\\
 S-751 06, Uppsala, Sweden
}
\end{center}

\vspace{0.4cm}

\begin{abstract}
In this paper, we  introduce   techniques for producing normal square-free monomial ideals from old
such ideals. These techniques are then used to investigate the normality of cover ideals under some graph operations. 
Square-free monomial ideals that come out as linear  combinations of two normal ideals are shown to be not necessarily normal; under such a  case we investigate the integral closedness of all powers of these ideals.  
\end{abstract}
\vspace{0.4cm}

\section{Introduction and preliminaries}
Let $R$ be a unitary commutative ring and $I$ an ideal in $R$. An element $f\in R$ is \emph{integral} over $I$, if there exists an equation 
 $$f^k+c_1f^{k-1}+\cdots +c_{k-1}f+c_k=0 ~~\mathrm{with} ~~ c_i\in I^i.$$
 The set of elements $\overline{I}$ in $R$ which are integral over $I$ is the 
 \emph{integral closure} of $I$. The ideal $I$ is 
 \emph{integrally closed}, if $I=\overline{I}$, and $I$ is 
 \emph{normal} if all powers of $I$ are integrally closed.
 This notion is linked to the graded algebras arising
from $I$ such as the Rees algebra $\mathrm{Rees}(I)=\oplus _{i\geq 0}I^{i}t^{i}$. It
is known that  if $I$ is  an ideal of a normal domain $R$, then $I$ is normal if and only if  $\mathrm{Rees}(I)$ is normal, see \cite[Theorem 4.3.17]{VI}. This brings up an importance of studying normality of ideals. 
When  $I$ is  a monomial ideal in a
polynomial ring $R$, then  $\overline{I}$ is the monomial ideal generated by all monomials $u \in R$ for which there exists an integer $k$ such that   $u^{k}\in I^{k}$, 
by  \cite[Theorem 1.4.2]{HH1}. In addition, it is well-known that every square-free monomial ideal  is  integrally closed, see \cite[Theorem 1.13]{ANR}.  
Appearing as edge
and cover ideals of graphs, the square-free monomial ideals play a key role in connecting commutative algebra and combinatorics, see \cite{EVY,MMV,VI3}.
The normality of such ideals has been of interest for many authors, see \cite%
{ANR,F,RV, SVV,VI,VI2}. For instance, it is shown in \cite{SVV} that the  edge ideals of bipartite graphs are normal. Also, it has been shown in \cite[%
Corollary 14.6.25]{VI} that cover ideals of perfect graphs are normal. Note that a graph is perfect if and only if it contains no odd cycle of length at least five, or its complement, as an induced subgraph, refer to \cite[Theorem 14.18]{BM}. 
In   \cite{ANR} it is shown that the cover ideals of odd cycles and wheel graphs
are normal. However, little is known for the cover ideals of imperfect  graphs. One goal of this paper is to investigate the normality of cover
ideals under graph operations. In particular, Theorem \ref{Th.Leaf}  says that if we
take any graph $G$ with normal cover ideal, and if we then add a leaf to $G$%
, then the cover ideal of the new graph is also normal. With this, and  Theorem 1.12 of \cite{ANR}, we prove that the cover ideal of a helm graph is
normal (Theorem \ref{Th. Helm}).

To  establish our results on cover ideals, we develop in Section
2 techniques that produce new normal square-free monomial
ideals from square-free monomial ideals already known to be
normal. In particular,  in Theorem \ref{Th.vI+wJ}  we prove
that the square-free monomial ideal $L:=x_{n}(IS:_{S}x_{n})+x_{n+1}I$ is normal in $S$%
, where $I\subset R=K[x_{1},\ldots ,x_{n}]$ is normal and $S=R[x_{n+1}]$.
Achieving the normality of $L$ has led the authors to investigate the  normality of square-free monomial ideals of the form $vI+wJ$. Initial
endeavours by the authors were to prove that such ideals are normal,  provided that $I$ and $J$ are normal. This assertion turned out to be not
true, see Remark \ref{linear}. Moreover, such linear combinations of normal ideals
came out to be a fruitful source of non-normal square-free monomial ideals.
Let $I_{n},J_{n}\subset R={\bQ}[x_{1},\ldots ,x_{n}]$ be the edge and  cover ideals of an $n$-cycle. Then $I_{n}$ is normal (see Corollary 6.1 and
Theorem 6.3 of \cite{SVV}) and also $J_{n}$ is normal (see \cite[Theorem 2.10%
]{VI4} and \cite[Theorem 1.10]{ANR}). In Section 4, we show that the ideal $%
L_{n}=x_{n+1}I_{n}+x_{n+2}J_{n}\subset R[x_{n+1},x_{n+2}]$ is not normal for
all $n\geq 4$. More specifically, we investigate the integral closedness of
every power $L_{n}^{i}$ for all $i\geq 2$.
In particular, we  prove, in  Theorems  
\ref{non-normal}  and \ref{normal}, that $L^2_5$ is not integrally closed but  $L^j_5$ is integrally closed for $j \geq 3$. We do  not  know any other cases of a monomial $J$ and integer $d$ such that $J^d$ is not integrally closed but $J^j$ is inegrally closed for $j>d$. The more typical behavior is illustrated in  Theorems  \ref{non-normal}  and 
\ref{non-normal-1}, where $L^j_n$ is not integrally closed $(j \geq 2$, $n=4$ or $n\geq 6$). 

Throughout this paper,  we denote the unique minimal set of monomial generators of a  monomial ideal $I$  by $\mathcal{G}(I)$. Also, in Sections 1-3, $R=K[x_1,\ldots, x_n]$ is a polynomial ring over a field $K$ and $x_1, \ldots, x_n$ are indeterminates. 
A simple graph $G$ means that $G$ has  no loop and no multiple edge. All graphs in this paper are undirected. 
Moreover, if $G$ is a finite simple graph, then $J(G)$  stands for the cover ideal of $G$. 

\section{Some results on the normality of monomial ideals}

In this section, we express several results on the normality of monomial ideals. To start with, recall from 
\cite[Definition 6.1.5]{VI} that if  $u=x_{1}^{a_{1}}\cdots x_{n}^{a_{n}}$ is  a
monomial in a polynomial ring $R=K[x_1, \ldots x_n]$ over a field $K$, then 
the \textit{support} of $u$ is given by $\mathrm{supp}(u):=\{x_{i}|~a_{i}>0%
\} $. 
   The following theorem is essential for us to show Corollary \ref{Cor_of_I+vh}. 
\begin{thm} \label{Th.I+vhR}
\label{I+vh}Let $I$ be a normal monomial ideal in $R$ and $h\in R$ a
monomial. Assume $v\in R$ is a square-free monomial with $\gcd (u,v)=1$ for
all $u\in \mathcal{G}(I)\cup \{h\}$. Then $L:=I+vhR$ is normal if and only
if $J:=I+hR$ is normal.
\end{thm}

\begin{proof}
Because $\gcd (u,v)=1$ for all $u\in \mathcal{G}(I)\cup \{h\}$, one may
assume that $v=x_{1}\cdots x_{m}\in K[x_{1},\ldots ,x_{m}]$ and $\mathcal{G}%
(I+hR):=\{u_{1},\ldots ,u_{s},h\}\subseteq K[x_{m+1},\ldots ,x_{n}]$ for
some positive integer $1\leq m<n$. The ideal $J$ is obtained from $L$\ by
making every variable $x_{i}$ (for $i=1,\ldots ,m$) equal  to $1$.
Proposition 12.2.3 of Villarreal \cite{VI} asserts that a normal monomial
ideal stays normal if we make any variable equals to $1$. This proves the
necessary part. Conversely, we prove that $\overline{L^{t}}=L^{t}$ for all
integers $t\geq 1$, provided that $J$ is normal. It suffices to show that $%
\overline{L^{t}}\subseteq L^{t}$. Let $\alpha$ be a monomial in $\overline{L^{t}}$ and
write $\alpha =v^{b}\delta $ with $v\nmid \delta $ and $\delta \in R$. By
\cite[Theorem 1.4.2]{HH1}, $\alpha ^{k}\in L^{tk}$ for some integer $k\geq 1$%
. Write%
\begin{equation}
v^{bk}\delta^{k}=\prod\limits_{i=1}^{s}u_{i}^{p_{i}}
v^{q+\varepsilon
}h^{q}\beta \text{,} \label{11}
\end{equation}%
with $\sum_{i=1}^{s}p_{i}+q=tk$, $\varepsilon \geq 0$, and $\beta $ is some
monomial in $R$ such that $v\nmid \beta $. 
Let $x_{l}\in 
\mathrm{supp}(v)\setminus \mathrm{supp}(\beta)$ for some $%
l\in \left\{ 1,\ldots ,m\right\} $. Now, if $x_{l}\in \mathrm{supp}(\delta)$ then, since $v$ is square-free, the $x_{l}$-degree in both sides of (%
\ref{11}) gives that $bk+jk=q+\varepsilon $, where $j=x_{l}$-$\deg (\delta )$%
. Hence, cancelling $v^{bk}$ from both sides of (\ref{11}) gives that $%
\delta ^{k}=\prod\nolimits_{i=1}^{s}u_{i}^{p_{i}}v^{q+\varepsilon
-bk}h^{q}\beta $. But $v$ is square-free; thus, $v$ must divide $\delta $, a
contradiction. Therefore, we proceed with the assumption that there exists $%
l\in \left\{ 1,\ldots ,m\right\} $ with $x_{l}\in 
\mathrm{supp}(v)\setminus 
\left(\mathrm{supp}(\beta)\cup 
\mathrm{supp}(\delta)\right) $. This
assumption, along with the hypothesis $\gcd (u,v)=1$ for all $u\in \mathcal{G%
}(I)\cup \{h\}$, forces  that $bk=q+\varepsilon $. Therefore, by (\ref{11}) we
obtain $\delta ^{k}=\prod\nolimits_{i=1}^{s}u_{i}^{p_{i}}h^{q}\beta $ $\in
J^{tk}$. This implies that $\delta $ $\in \overline{J^{t}}$. Since $J$ is
normal, then $\overline{J^{t}}=J^{t}$, thus, $\delta $ $\in J^{t}$; hence, 
write%
\begin{equation}
\delta =\prod\limits_{i=1}^{s}u_{i}^{l_{i}}h^{z}\gamma \text{,} \label{22}
\end{equation}%
with $\sum_{i=1}^{s}l_{i}+z=t$, $z\geq 0$, and $\gamma $ is some monomial in 
$R$. Note $v\nmid \gamma $ as $v\nmid \delta $. Now, since $v^{bk}\delta
^{k}\in L^{tk}$ then by (\ref{22}) we get $\prod%
\nolimits_{i=1}^{s}u_{i}^{l_{i}k}v^{bk}h^{zk}\gamma ^{k}\in L^{tk}=\left(
I+vhR\right) ^{tk}$. Therefore, we conclude that $bk\geq zk$, that is, $%
b\geq z$. Thus, $v^{b}\delta
=\prod\nolimits_{i=1}^{s}u_{i}^{l_{i}}v^{b}h^{z}\gamma \in \left(
I+vhR\right) ^{t}$, and the proof is done.
\end{proof}


Setting $h=1$ in Theorem \ref{Th.I+vhR}, we obtain the following corollary.

\begin{cor} \label{Cor_of_I+vh} 
Let $I$ be a normal monomial ideal in $R$ and $v\in R$ a
square-free monomial with $\mathrm{gcd}(u,v)=1$ for all $u\in \mathcal{G}(I)$%
. Then $I+vR$ is normal.
\end{cor}

Inductively and in view of  \cite[Proposition 12.2.3]{VI}, the following proposition follows directly from 
Corollary \ref{Cor_of_I+vh}.

\begin{pro} \label{Pro.I+J}
\label{Prop_of_I+vh} 
Any square-free monomial ideal of $R$ with a set of pairwise relatively prime generators is normal. Assume $%
\mathcal{G}(I)\subset K[x_{1},\ldots ,x_{m}]$ and $J$ is square-free with $\mathcal{G}%
(J)\subset K[x_{m+1},\ldots ,x_{n}]$ with $%
1\leq m<n$ and the generators of $J$ are pairwise relatively prime. Then $%
I+J $ is normal if and only if $I$ is  normal. 
\end{pro}


\begin{rem}
It should be noted that Proposition \ref{Pro.I+J} may be false if we drop the condition that the generators of $J$ are pairwise relatively prime. To see a counterexample, let  $I:=(x_1x_2x_4, x_1x_3x_4, x_1x_3x_5, x_2x_3x_5, x_2x_4x_5) \subset K[x_1, \ldots, x_5]$ and $J:=(y_1y_2y_4, y_1y_3y_4, y_1y_3y_5, y_2y_3y_5, y_2y_4y_5) \subset  K[y_1, \ldots, y_5]$ be  the cover ideals of the odd cycle graphs  $G$ and $H$,  respectively, where  
$V(G)=\{x_1,x_2,x_3,x_4,x_5\}$   with 
$E(G)=\{\{x_i, x_{i+1}\}\}_{i=1}^5,$  and  
$V(H)=\{y_1,y_2,y_3,y_4,y_5\}$ with 
$E(H)=\{\{y_i, y_{i+1}\}\}_{i=1}^5$ and $x_6$ represents $x_1$, and $y_6$ represents $y_1$. 
It follows from \cite[Theorem 1.10]{ANR} that $I$ and $J$ are normal monomial ideals. Now, set $\alpha:=\displaystyle\prod_{i=1}^5x_i\prod_{i=1}^5y_i$ in $$Q:=IK[x_1, \ldots, x_5,y_1, \ldots, y_5] +JK[x_1, \ldots, x_5, y_1, \ldots, y_5].$$ Because 
$$\alpha^{2}=(x_1x_2x_4)(x_1x_3x_5)(x_2x_3x_5)x_4
(y_1y_2y_4)(y_1y_3y_5)(y_2y_3y_5)y_4\in Q^6,$$ 
one can deduce that $\alpha \in \overline{Q^3}$. On the other hand, it is easy to check that $\alpha \notin Q^3$. This means that $Q$ is non-normal. 

\end{rem}


Assume a graph $H$ is obtained from a graph $G$ by connecting all the vertices of $G$ with a new vertex. In \cite[Theorem 1.6]{ANR} it is proved
that the cover ideal of the graph $H$ is normal provided that the cover ideal of the  graph $G$ is normal. The proof relies  on \cite[Theorem 1.4]{ANR} in which it is proved that $vI+hR$ is normal provided that $I$ is a normal
monomial ideal, $h\in I$, and $v$ is a square-free monomial which is relatively prime to every generator of $I$. In 
Theorem \ref{New.Th.ANR}, we give a
generalization   of  \cite[Theorem 1.4]{ANR}.  The following 
lemma is needed in the proofs of  Theorems    \ref{New.Th.ANR} and  \ref{Th.vI+wJ}. It is a straightforward  application of  \cite[Theorem  1.4.2]{HH1}.   

\begin{lem}
\label{Lem.Intersection} 
Suppose that $I$ and $J$ are two normal monomial
 deals in $R$ such that   $\mathrm{gcd}(u,v)=1$ for all 
$u\in \mathcal{G}(I)$ and $v\in \mathcal{G}(J)$. Then $I\cap J=IJ$ is normal.
\end{lem}


\begin{thm} \label{New.Th.ANR} Let $I$ be a normal monomial ideal of $R$, and $J$ be a normal square-free monomial ideal of $R$   whose generators are pairwise relatively prime. Let $%
h\in I$ a monomial and suppose $\mathrm{gcd}(u,v)=1$ for all $u\in \mathcal{G%
}(I)\cup \{h\}$ and $v\in \mathcal{G}(J)$. Then $L:=JI+hR$ is normal.
\end{thm}

\begin{proof}
Assume that $\mathcal{G}(I):=\{u_{1},\ldots ,u_{r}\}$ and $\mathcal{G}%
(J):=\{v_{1},\ldots ,v_{s}\}$. As before, it suffices to show that $%
\overline{L^{t}}\subseteq L^{t}$ for all $t\geq 1$.
 Let $\alpha$ be a monomial in $\overline{L^{t}}$. Then we have $\alpha ^{k}\in
(JI+hR)^{tk}$ for some positive integer $k$. The binomial expansion implies
that $\alpha ^{k}\in (JI)^{p}(hR)^{q}$ for some nonnegative integers $p$ and 
$q$ with $p+q=tk$. Choose $q$ to be the minimal according to this
membership. If $q=0$, then $p=tk$; hence,  $\alpha ^{k}\in (JI)^{tk}$ and $\alpha \in \overline{(JI)^{t}}$. In the light of $I$ and $J$
being  normal and $\mathrm{gcd}(u,v)=1$ for all $u\in \mathcal{G}(I)$ and $v\in 
\mathcal{G}(J)$, Lemma \ref{Lem.Intersection} yields that $JI$ is normal.
and so $\alpha \in (JI)^{t}$. Thus, we have $\alpha \in L^{t}$ and the proof
is finished. Therefore, let $q\geq 1$. Since $\alpha ^{k}\in (JI)^{p}(hR)^{q}$,
we get the following equality 
\begin{equation}
\alpha ^{k}=v_{1}^{z_{1}}\cdots v_{s}^{z_{s}}u_{1}^{p_{1}}\cdots
u_{r}^{p_{r}}h^{q}\beta , \label{Eq.3}
\end{equation}%
with $\sum_{i=1}^{s}z_{i}=p\sum_{i=1}^{r}p_{i}$, and $\beta $  some
monomial in $R$. Since  $h\in I$, then $h=h^{\prime }u_{\lambda }$ for some 
$1\leq \lambda \leq r$ and monomial $h^{\prime }$ in $R$. If $\beta \in J$,
then $\beta =v_{\theta }\beta ^{\prime }$ for some $1\leq \theta \leq s$ and
monomial $\beta ^{\prime }$ in $R$. Hence, $(\ref{Eq.3})$ can be rewritten
as follows 
\begin{equation*}
\alpha ^{k}=v_{1}^{z_{1}}\cdots v_{s}^{z_{s}}v_{\theta }u_{1}^{p_{1}}\cdots
u_{r}^{p_{r}}u_{\lambda }h^{q-1}h^{\prime }\beta ^{\prime }.
\end{equation*}%
This leads to $\alpha ^{k}\in (JI)^{p+1}(hR)^{q-1}$, which contradicts the
minimality of $q$. Accordingly, one can assume in $(\ref{Eq.3})$ that $%
v_{i}\nmid \beta $ for each $i=1,\ldots ,s$. Let $x_{j_{i}}|v_{i}$ but $%
x_{j_{i}}\nmid \beta $ for each $i=1,\ldots ,s$. Note that $x_{j_{1}},\ldots
,x_{j_{s}}$ are distinct since $\mathrm{gcd}(v_{i},v_{j})=1$ for any $1\leq
i\neq j\leq s$. Write $\alpha =x_{j_{1}}^{b_{1}}\cdots
x_{j_{s}}^{b_{s}}\delta $ with $x_{j_{d}}\nmid \delta $ for each $d=1,\ldots
,s$. This gives rise to the following equality 
\begin{equation}
x_{j_{1}}^{b_{1}k}\cdots x_{j_{s}}^{b_{s}k}\delta ^{k}=v_{1}^{z_{1}}\cdots
v_{s}^{z_{s}}u_{1}^{p_{1}}\cdots u_{r}^{p_{r}}h^{q}\beta . \label{Eq.4}
\end{equation}%
In view of $(\ref{Eq.4})$, we have $b_{i}k=z_{i}$ for each $i=1,\ldots ,s$.
Set $b:=b_{1}+\cdots +b_{s}$. Thanks to $\sum_{i=1}^{s}z_{i}=p$, one has $%
bk=p$, and so $q=(t-b)k$. We thus have $\alpha ^{k}\in
(JI)^{bk}(hR)^{(t-b)k} $. Consequently, $\alpha \in \overline{%
(JI)^{b}(hR)^{t-b}}$. It follows from Lemma \ref{Lem.Intersection} that $%
\overline{(JI)^{b}(hR)^{t-b}}=(JI)^{b}(hR)^{t-b}$. This yields that $\alpha
\in (JI)^{b}(hR)^{t-b}\subseteq (JI+hR)^{t}=L^{t}$, as required.
\end{proof}


The subsequent theorem is one of the main results in this paper which is used in   proving Theorem \ref{Th.Leaf}. 
\begin{thm}\label{Th.vI+wJ} 
Let $I$ be a normal square-free monomial ideal in $%
R=K[x_{1},\ldots ,x_{n}]$ with $\mathcal{G}(I) \subset R$. Then the ideal $L:=IS\cap (x_{n},x_{n+1})\subset  S=R[x_{n+1}]$ is normal.
\end{thm}

\begin{proof}
Since $I\cap (x_{n})=x_{n}\left( I:_{S}x_{n}\right) $ and $I\cap
(x_{n+1})=x_{n+1}I$, one can conclude that $L=x_{n}\left( I:_{S}x_{n}\right)
+x_{n+1}I$. To simplify notation, set $F:=\left( I:_{S}x_{n}\right) $ and $%
L=x_{n}F+x_{n+1}I$. Since $I$\ is square-free, then the ideal $F$ is
obtained from $I$\ by making the variable $x_{n}$ equal  to $1$. Therefore, $%
F$ is normal by  virtue of Proposition 12.2.3 in  \cite{VI},
and hence $x_{n}F$ is also normal. Let $\mathcal{G}(I)=\{g_{1},\ldots
,g_{d},g_{d+1},\ldots ,g_{s}\}$ with $x_{n}\mid g_{j}$ for $j=1,\ldots ,d$
and $x_{n}\nmid g_{j}$ for $j=d+1,\ldots ,s$. Then a generating set (not
necessarily minimal) of $F$ is given by $\{f_{1},\ldots
,f_{d},f_{d+1},\ldots ,f_{s}\}$ with $f_{i}x_{n}=g_{i}$ for $i=1,\ldots ,d$
and $f_{i}=g_{i}$ for $i=d+1,\ldots ,s$. Note that $x_{n}F+I=I$; therefore, $%
x_{n}F+I$ is normal, that is, $\overline{\left( x_{n}F+I\right) ^{t}}=\left(
x_{n}F+I\right) ^{t}$ for all $t\geq 1$. Our goal is to show that $\overline{%
L^{t}}=L^{t}$ for all $t\geq 1$. 

Let $\alpha$ be a monomial in $\overline{L^{t}}$ and write 
$\alpha =x_{n+1}^{b}\delta $ for some integer $b$ and some monomial $\delta
\in R$ with $x_{n+1}\nmid \delta $. As $\alpha \in \overline{L^{t}}$,  \cite[Theorem 1.4.2]{HH1} implies that 
 $%
\alpha ^{k}\in L^{tk}=\left( x_{n}F+x_{n+1}I\right) ^{tk}$ for some integer $%
k$; therefore, $\alpha ^{k}\in \left( x_{n}F\right) ^{p}\left(
x_{n+1}I\right) ^{q}$ for some integers $p$ and $q$ with $p+q=tk$. Assume $q$
is maximal according to this membership. Note that if $p=0$,  then $\alpha
^{k}\in \left( x_{n+1}I\right) ^{tk}$, and hence $\alpha \in \overline{\left(
x_{n+1}I\right) ^{t}}=\left( x_{n+1}I\right) ^{t}\subset L^{t}$. Henceforth,
assume $p>0$. Similarly, and since $x_{n}F$ is normal, we may also assume $%
q>0$. Write%
\begin{equation}
\alpha ^{k}=x_{n+1}^{bk}\delta ^{k}=\prod\limits_{i=1}^{s}f_{i}^{p_{i}}%
\text{ }x_{n}^{p}\text{ }\prod\limits_{j=1}^{s}g_{j}^{q_{j}}\text{ }%
x_{n+1}^{q}\text{ }\beta \text{,}  \label{111}
\end{equation}%
with $\sum\nolimits_{i=1}^{s}p_{i}=p$, $\sum\nolimits_{j=1}^{s}q_{j}=q$,
and $\beta $   some monomial in $S$. If $x_{n+1}\mid \beta $, then we get a
contradiction to the maximality of $q$ since either $f_{i}x_{n}=g_{i}$ or $%
f_{i}=g_{i}$. Therefore, we may assume in (\ref{111})\ that $x_{n+1}\nmid
\beta $, and thus we can conclude that $q=bk$ and also that%
\begin{equation}
\delta ^{k}=\prod\limits_{i=1}^{s}f_{i}^{p_{i}}\text{ }x_{n}^{p}\text{ }%
\prod\limits_{j=1}^{s}g_{j}^{q_{j}}\text{ }\beta \in \left( x_{n}F+I\right)
^{tk}\text{.}  \label{222}
\end{equation}%
Therefore, $\delta \in $ $\overline{\left( x_{n}F+I\right) ^{t}}=\left( x_{n}F+I\right) ^{t}$. Thus, let $\delta \in \left(
x_{n}F\right) ^{l}I^{h}$ with $l+h=t$ and $l$ being maximal with respect to this membership. Note that if $h=0$, then $\delta \in \left( x_{n}F\right)
^{t}\subset L^{t}$. Henceforth, assume $h>0$. Write%
$$
\delta =\prod\limits_{i=1}^{s}f_{i}^{l_{i}}\text{ }x_{n}^{l}\text{ }%
\prod\limits_{j=1}^{s}g_{j}^{h_{j}}\text{ }\gamma \text{,}  
$$
with $\sum\nolimits_{i=1}^{s}l_{i}=l$, $\sum\nolimits_{j=1}^{s}h_{j}=h$,
and $\gamma $ is some monomial in $S$. Note that if $d=0$, then $L=I\cap (x_n, x_{n+1})$ is normal by Lemma \ref{Lem.Intersection}. In addition, note that if $d=s$, 
then  $L=I$, which is normal. 
Henceforth, assume that $s>d>0$. 
 The maximality of $l$ implies that $%
h_{j}=0$ for $j=1,\ldots ,d$ and also that $x_{n}\nmid \gamma $ since $%
g_{j}=f_{j}$ for $j=d+1,\ldots ,s$. Hence, $x_{n}$-$\deg (\delta )=l$. Since 
$p=(t-b)k$, then (\ref{222}) gives that $x_{n}$-$\deg (\delta )\geq t-b$, that
is, $l\geq t-b$; thus, $b\geq h$. Therefore, since $\delta \in \left(
x_{n}F\right) ^{l}I^{h}$, one can deduce that  $$x_{n+1}^{b}\delta \in \left(
x_{n}F\right) ^{l}\left( x_{n+1}I\right) ^{h}\subseteq \left(
x_{n}F+x_{n+1}I\right) ^{l+h}=L^{t},$$ which finishes the proof.
\end{proof}


\begin{rem} \label{linear}
As the reader may notice, in the proof of Theorem \ref{Th.vI+wJ}
we proved that the linear combination $x_{n}F+x_{n+1}I$ is normal, where $%
I\subseteq F$. Initial endeavours of the authors were to prove a more
general result, that is, investigating the normality of square-free monomial ideals
resulting from linear combinations $vF+wI$, where $F$ and $I$ are normal
ideals, $v$ and $w$ are square-free monomials with $\mathrm{gcd}(v,f)=1$ for all $f\in \mathcal{G}(F)$ and 
$\mathrm{gcd}(w,g)=1$ for all $g\in \mathcal{G}(I)$. Since $vF+wI$ is square-free, then it is integrally closed. However,
we found that one cannot guarantee the normality of $vF+wI$, even though one
has $I\subseteq F$ or $\mathcal{G}(I)\subset \mathcal{G}(F)$. In this remark
we demonstrate an example supporting this assertion. Let $%
F=(x_{1}x_{2}x_{4},x_{1}x_{3}x_{5},x_{2}x_{3},x_{2}x_{5},x_{3}x_{4})$ and $%
I=(x_{1}x_{2}x_{4},x_{1}x_{3}x_{5})$ in $R=K[x_{1},\ldots ,x_{7}]$, and let $%
v=x_{6}$ and $w=x_{7}$. Then $T:=vF+wI=$
\[
=(x_{1}x_{2}x_{4}x_{6},x_{1}x_{3}x_{5}x_{6},x_{2}x_{3}x_{6},x_{2}x_{5}x_{6},x_{3}x_{4}x_{6},x_{1}x_{2}x_{4}x_{7},x_{1}x_{3}x_{5}x_{7}).
\]%
The normality of $I$ can be deduced from Proposition \ref{Pro.I+J} and \ref{Lem.Intersection}.  In addition, using
\emph{Normaliz} \cite{N} yields that $F$ is a normal monomial ideal. Now, put $\alpha
:=x_{1}x_{2}x_{3}x_{4}x_{5}x_{6}x_{7}$. Direct computations show that $%
\alpha \notin T^{2}$. Since%
\[
\alpha ^{2}=\left( x_{2}x_{5}x_{6}\right) \left( x_{3}x_{4}x_{6}\right)
\left( x_{1}x_{2}x_{4}x_{7}\right) \left( x_{1}x_{3}x_{5}x_{7}\right) \in
T^{4}\text{,}
\]%
we conclude that  $\alpha \in \overline{T^{2}}\backslash T^{2}$, that is, $%
vF+wI$ is not normal.
\end{rem}

\begin{cor} 
Let $I$ be a normal square-free monomial ideal in $%
R=K[x_{1},\ldots ,x_{n}]$ with $\mathcal{G}(I) \subset R$. Then the ideal $L:=IS\cap (x_{n},x_{n+1}\cdots x_m)\subset S=R[x_{n+1}, \ldots, x_m]$ is normal.
\end{cor}
\begin{proof}
Since  $I\cap (x_{n},x_{n+1}\cdots x_m)=I\cap (x_n, x_{n+1})\cap (x_n, x_{n+2}) \cap \cdots \cap (x_n, x_{m}),$ this claim follows at once from  Theorem \ref{Th.vI+wJ}.
\end{proof}


The above corollary motivates for the following two questions.

\begin{question}\label{Question1}
Let $I$ be a normal square-free monomial ideal in $%
R=K[x_{1},\ldots ,x_{n}]$  with $\mathcal{G}(I) \subset R$, and $\{i_1, \ldots, i_r\} \subseteq \{1, \ldots, n\}$ with $r>1$. Then, in general, can we deduce that 
 the ideal $IS\cap (x_{i_1}\ldots x_{i_r}, x_{n+1})\subset S=R[x_{n+1}]$ is normal?
\end{question}

\begin{question}\label{Question2}
Let $I$ be a normal square-free monomial ideal in $R$. Then, in general, can one conclude that   $IS\cap (x_{n},x^{\ell}_{n+1})\subset S=R[x_{n+1}]$, with $\ell >1$, is normal?
\end{question}

We show that Question \ref{Question1} has a negative answer, while we leave Question \ref{Question2} open.    
For this purpose, we provide  a counterexample. Firstly, one should recall that, given a graph $G$, 
if $v$ is a vertex of $G$, we may obtain a graph on $n-1$
vertices by deleting from $G$ the vertex $v$ together with all the edges incident with $v$. The resulting graph is denoted by $G\setminus v$.  In the sequel, 
 we consider the graph   \cite[Section 2]{KSS}, which is described in the following way. For a positive integer $n$, let $[n]$ denote the set $\{0, \ldots, n-1\}$. Assume that  $P_n$ is  a path with vertex set $[n]$, with vertices in the increasing order along $P_n$. Let also   $K_3$ be  the complete graph whose vertex set is the group ${\bZ}_3$. For $n\geq 4$, we define $H_n$ as the graph obtained from the Cartesian product $P_n\Box K_3$ by adding the three edges joining $(0,j)$ to $(n-1,-j)$ for $j\in {\bZ}_3$. Figure  below is the graph of $H_4$ in 
\cite{KSS}. 

\begin{figure}[h!]
\centering
\begin{tikzpicture}
\foreach \x/\xtext  in {0,1,2,3}{
\foreach \y/\ytext  in {0/2,1/0,2/1}{
\filldraw (\x cm,\y cm) circle [radius=.08cm] ;
\draw (0,\y) to (1,\y) to (2,\y) to (3,\y);
\draw (\x,0) to (\x, 1) to (\x,2);
\draw (\x, 0) to [in=-120 , out=120] (\x, 2);
\pgfmathsetmacro{\xvar}{\x + .012cm};
\pgfmathsetmacro{\yvar}{\y + .005cm};
\draw (\xvar, \yvar) node {{\scriptsize{$v_{\xtext,\ytext}$}}};
}}
\draw (0,2)  to[in=130 , out=-130 ] (-.2,-.2) to [in=-150 , out=-45 ] (3,0);
\draw (3,2)  to[in=50 , out=-50 ] (3.2,-.2) to [in=-30 , out=-135 ] (0,0);
\draw (0,1) to [in=-150 , out=-30] (3,1);
\end{tikzpicture}
\caption{$H_4$}
\end{figure}

Assume that  $J(H_4)$ denotes  the cover ideal of  $H_4$ in the polynomial ring $S = K[x_\alpha : \alpha \in V(H_4)]$ over a field $K$.  Now,  put $G:=H_4\setminus v_{0,1}$. Let $J(G)$ denote the cover ideal of the graph $G$. 
It should be noted that 
\begin{align*}
J(H_4)=& J(G) \cap (x_{v_{0,1}}, x_{v_{1,1}}) \cap (x_{v_{0,1}}, x_{v_{0,0}}) \cap ((x_{v_{0,1}}, x_{v_{0,2}})\cap (x_{v_{0,1}}, x_{v_{3,2}})\\
=& J(G) \cap (x_{v_{0,1}}, x_{v_{1,1}}x_{v_{0,0}}x_{v_{0,2}}x_{v_{3,2}}).
\end{align*}

It follows from Normaliz \cite{N} that $J(G)$ is a normal monomial ideal.
 While, according to \cite[Page 21]{ANR},  $J(H_4)$ is non-normal.

\section{An application to the cover ideals of imperfect graphs}

The aim of this
section is to apply some of the results of the previous section to
cover ideals of imperfect graphs. To achieve this, we need to recall some definitions from graph theory. \par 

Let $G$ be a finite simple graph with the vertex set $V(G)$ and 
the edge set $E(G)$. A subset $W\subseteq V(G)$ is called a 
\emph{vertex cover} of $G$ if it intersects any edge of $G$. Furthermore,  $W$ is called a \emph{minimal vertex cover} of $G$ if it is a vertex cover and no proper subset of $W$ is a vertex cover of $G$. Let $W_1, \ldots, W_r$ be the minimal vertex covers of the graph $G$. Then, the \emph{cover ideal} of $G$, denoted by $J(G)$, is defined as $J(G)=(X_{W_1}, \ldots, X_{W_r}),$ where 
$X_{W_j}=\prod_{t\in W_j}x_t$ for each $j=1, \ldots, r$.
 For more information about cover ideals and the other kinds of cover sets see \cite{HH1, KHN1, N4, SNK}. 
 
\begin{defn}
Let $G=(V(G), E(G))$  be  a  finite simple graph. A $k$-coloring of $G$ is any  partition of $V(G)=C_1 \cup \cdots \cup C_k$ into $k$ disjoint sets such that for each $e\in E(G)$, one has $e\nsubseteq C_i$ for all $i=1, \ldots, k$. The chromatic number of $G$, denoted by $\chi(G)$, is the minimal $k$ such that $G$ has a $k$-coloring. 
\end{defn}
\begin{defn}
Let $G=(V(G), E(G))$  be  a  finite simple graph. The graph $G$ is called perfect if for any induced subgraph $G_S$, with $S\subseteq V(G)$, we have $\chi(G_S)=\omega(G_S)$, where 
$\omega(H)$ denotes the size of the largest clique of a graph $H$.
\end{defn}
It is well-known that  a graph is perfect if and only if it contains no odd cycle of length at least five, or its complement, as an induced subgraph, by  \cite[Theorem 14.18]{BM}.
Moreover, it has been shown in \cite[%
Corollary 14.6.25]{VI} that cover ideals of perfect graphs are normal.  However, little is known about the normality of cover ideals of imperfect graphs.

\par First we prove a result on adding a leaf to a graph with normal cover ideal.  In this section, $[n]=\{1, \ldots, n\}$.

\begin{thm}\label{Th.Leaf}
Let $G=(V(G), E(G))$ and  $H=(V(H), E(H))$ be  finite simple graphs such that $V(H)=V(G)\cup \{w\}$ with $w\notin V(G)$,  and $E(H)=E(G) \cup \{\{v,w\}\}$ for some vertex $v\in V(G)$. Let $J(G)$ and  $J(H)$ be  the cover ideals of the graphs $G$ and   $H$, respectively.  If $J(G)$ is normal,  then $J(H)$ is  normal.
\end{thm}
\begin{proof}
Suppose that $J(G)$ is normal. Without loss of generality, one may  assume  that  $V(G)=[n]$, $V(H)=V(G)\cup \{n+1\}$, and  $E(H)=E(G)\cup \{\{n, n+1\}\}$.   Since $J(H)=J(G) \cap (x_n, x_{n+1})$, the claim  is a straightforward consequence of Theorem \ref{Th.vI+wJ}, that is, $J(H)$ is normal. 
\end{proof}

Here, we want to explore the normality of the cover ideals of helm graphs $H_n$ for $n\geq 5$ odd. For this purpose, one requires to recall  the following definitions and a theorem.

\begin{defn} 
A wheel graph $W_n$ of order $n$ is a graph that contains a cycle of order $n-1$, and for  which every vertex in the cycle is connected to one other vertex which is known as the hub. The edges of a wheel which include the hub are called spokes.
\end{defn}

\begin{thm}  \cite[Theorem 1.12]{ANR}\label{Th. Wheel}
Suppose that $W_{2n}$ is  a wheel graph  of order $2n$ on the  vertex set $[2n]$. Then $J(W_{2n})$  is normal. 
\end{thm}

\begin{defn}
The helm graph $H_n$, which  has $2n+1$ vertices, is the graph obtained from a wheel graph $W_{n+1}$ of order $n+1$   by adjoining a pendant edge at each node of the outer circuit of the  wheel graph $W_{n+1}$.  
\end{defn}
As an application of  Theorem \ref{Th.Leaf}, we illustrate that   every cover ideal of  helm graphs $H_n$ for $n\geq 5$  odd,  is normal. It should be noted that based on \cite[Theorem 14.18]{BM},  a graph is perfect if and only if it contains no odd cycle of length at least five, or its complement, as an induced subgraph. Since  $H_{2n+1}$  with $n\geq 2$, contains an induced odd cycle of length $n\geq 5$, this graph is imperfect, and so the normality of its cover ideal is of special interest.

Also, recall that if  $I$ is  an ideal in a commutative Noetherian ring $S$, then  $I$ is said to have the \emph{persistence property} if $\mathrm{Ass}_S(S/I^k)\subseteq \mathrm{Ass}_S(S/I^{k+1})$ 
for all positive integers $k$. Moreover, an ideal $I$ satisfies the \emph{strong persistence property} if $(I^{k+1}:_S I)=I^k$ for all positive integers $k$. Specially, it is well-known that the strong persistence property implies  the persistence property, see \cite{HQ}.
\begin{thm}\label{Th. Helm}
Suppose that $H_{2n+1}$  with $n\geq 2$ is  a helm  graph on the  vertex set $[4n+3]$. Then $J(H_{2n+1})$  is normal. 
 Therefore, it has the strong persistence property, and hence the persistence property.
\end{thm}
\begin{proof}
Label the vertices of  $C_{2n+1}$ by $1,\ldots, 2n+1$  in counterclockwise order, and the hub by $2n+2$, as shown in figure below,   such that we have 
$$E(H_{2n+1})=E(W_{2n+2})\cup \{\{i, 2n+2+i\}~:~ i=1, \ldots, 2n+1\}.$$
$\hspace{1cm}$

\begin{figure}[h!]

\scalebox{.8}  
{
\begin{pspicture}(-1,-2.558125)(5.6028123,2.558125)
\psdots[dotsize=0.12](2.7609375,1.0596875)
\psdots[dotsize=0.12](3.9409375,0.4596875)
\psdots[dotsize=0.12](1.5609375,0.4396875)
\psdots[dotsize=0.12](2.7409375,0.0396875)
\psdots[dotsize=0.12](3.3409376,-0.7403125)
\psdots[dotsize=0.12](2.1409376,-0.7403125)
\psline[linewidth=0.024cm](2.7409375,1.0796875)(1.5409375,0.4396875)
\psline[linewidth=0.024cm](2.7409375,1.0796875)(3.9209375,0.4796875)
\psline[linewidth=0.024cm](3.9209375,0.4796875)(3.3409376,-0.7203125)
\psline[linewidth=0.024cm](2.7209375,0.0996875)(2.7009375,0.0396875)
\psline[linewidth=0.024cm](2.7409375,0.0596875)(2.7409375,1.0796875)
\psline[linewidth=0.024cm](2.7209375,0.0996875)(2.7609375,0.0996875)
\psline[linewidth=0.024cm](2.7209375,0.0596875)(3.9409375,0.4796875)
\psline[linewidth=0.024cm](2.7209375,0.0596875)(1.5809375,0.4396875)
\psline[linewidth=0.024cm](1.5809375,0.4396875)(1.5409375,0.4396875)
\psline[linewidth=0.024cm](1.5409375,0.4796875)(2.1209376,-0.7203125)
\psline[linewidth=0.024cm](2.1209376,-0.7203125)(2.7409375,0.0396875)
\psline[linewidth=0.024cm](2.7409375,0.0596875)(3.3209374,-0.7003125)
\psdots[dotsize=0.12](4.9609375,0.8396875)
\psdots[dotsize=0.12](4.1209373,-1.5403125)
\psdots[dotsize=0.12](1.3409375,-1.5203125)
\psdots[dotsize=0.12](0.5409375,0.8596875)
\psdots[dotsize=0.12](2.7409375,2.0196874)
\psline[linewidth=0.024cm](2.7409375,2.0396874)(2.7409375,1.0796875)
\psline[linewidth=0.024cm](3.9209375,0.4796875)(4.9609375,0.8396875)
\psline[linewidth=0.024cm](4.1209373,-1.5403125)(3.3209374,-0.6803125)
\psline[linewidth=0.024cm](1.3209375,-1.4803125)(1.3209375,-1.5203125)
\psline[linewidth=0.024cm](1.3409375,-1.5003124)(2.1209376,-0.7403125)
\psline[linewidth=0.024cm](0.5209375,0.8796875)(1.5209374,0.4596875)
\psline[linewidth=0.024cm,linestyle=dashed,dash=0.16cm 0.16cm](2.1209376,-0.7203125)(3.3009374,-0.7203125)
\usefont{T1}{ptm}{m}{n}
\rput(3.0323439,1.2696875){$1$}
\usefont{T1}{ptm}{m}{n}
\rput(2.7923439,2.3696876){$2n+3$}
\usefont{T1}{ptm}{m}{n}
\rput(1.5323437,0.7696875){$2$}
\usefont{T1}{ptm}{m}{n}
\rput(0.53234375,1.1496875){$2n+4$}
\usefont{T1}{ptm}{m}{n}
\rput(4.912344,1.2096875){$4n+3$}
\usefont{T1}{ptm}{m}{n}
\rput(3.1723437,0.5096875){ \small\emph{2n+2}}
\usefont{T1}{ptm}{m}{n}
\rput(3.9223437,-0.6903125){$2n$}
\usefont{T1}{ptm}{m}{n}
\rput(4.912344,-1.4303125){$4n+2$}
\usefont{T1}{ptm}{m}{n}
\rput(1.6123438,-0.7103125){$3$}
\usefont{T1}{ptm}{m}{n}
\rput(0.51234376,-1.4503125){$2n+5$}
\usefont{T1}{ptm}{m}{n}
\usefont{T1}{ptm}{m}{n}
\rput(4.472344,0.2496875){$2n+1$}
\end{pspicture} 
}
\caption{$H_{2n+1}$}
\end{figure}

We can now combine together  Theorem \ref{Th. Wheel}  and    the iteration of   Theorem \ref{Th.Leaf}  to obtain the 
normality of $J(H_{2n+1})$. The last  assertion follows readily from the normality of $J(H_{2n+1})$.
\end{proof}


\section{An argument on the normality of linear combinations of two normal ideals}

 Let $R={\bQ}[x_1,\ldots,x_{n+2}]$, let $I$ and $J$ be two square-free monomial ideals in ${\bQ}[x_1,\ldots,x_n]$, and let $L=x_{n+1}IR+x_{n+2}JR$.
As we have mentioned in Remark \ref{linear}, this construction is a fruitful source of interesting square-free monomial ideals $L$ which are not normal.

In this section, we  investigate the case where $I$ is the edge ideal
$$I_n=(x_1x_2,x_2x_3,\ldots, x_{n-1}x_n,x_nx_1),$$ of an $n$-cycle $C_n, n\geq 3$ and  $J$ is the cover ideal $J_n=(x_n,x_1)\cap_{1\leq i\leq n-1}(x_i,x_{i+1})$  of $I_n$.
We  will write $L_n$ instead of $L$,
so that $L_n=x_{n+1}I_nR+x_{n+2}J_nR$. By abuse of notation, we  will write more simply $L_n=x_{n+1}I_n+x_{n+2}J_n$. Also all statements about $L_n, I_n, J_n$ will take place in the
ring $R={\bQ}[x_1,\ldots,x_{n+2}]$, which for simplicity of notation, we  might not explicitly mention. In a similar vein  $\overline{L_n^i}\backslash L_n^i$ literally means the complement of $L_n^i$ in
$\overline{L_n^i}$, i.e. those elements of $\overline{L_n^i}$ not in $L_n^i$. But  computationally we like to think of it as all the monomials in
$\overline{L_n^i}$ that are not in $L_n^i$. Macaulay2 calculation shows that $L_3$ is normal. In the sequel,  we will
always have $n\geq 4$. 

Our  first  conclusion is that  $L_n$ is not normal
for $n\geq 4$. More specifically, $L_n^2$ is not integrally closed for any $n\geq 4$ (Theorem \ref{non-normal}). If $i\geq 3$,  then $L_n^i$ is not integrally closed for even $n\geq 4$, and for odd $n\geq 7$
(Theorem \ref{non-normal-1}).
If $n=5$,  then $L_5^i$ is integrally closed for $i\geq 3$ (Theorem \ref{normal}). We  used the Hilbert basis of the Rees cone to prove the last result. Thus, we  characterize  all cases  when $L_n^i$ is integrally closed, i.e., when $\overline{L^i_n}/L^i_n=0$.

We finally show that  $\overline{L_n^2}/L_n^2$ is not a finite dimensional vector space for all odd $n$, $n\geq 5$ (Theorem \ref{non-normal-3}). 

$\vspace{.03cm}$
\subsection{The case of $L_n^2$} 
$\vspace{.1cm}$

In this subsection, we  prove

\begin{thm}\label{non-normal} If $n\geq 4$, then $L_n^2$ is not integrally closed.
\end{thm}

\begin{proof}
Let $f=x_1x_2\cdots x_nx_{n+1}x_{n+2}$. Then  we  claim
    that $f\notin L_n^2$. If we assume the opposite, then $f$ would be the product of one minimal generator of $x_{n+1}I_n$ and one minimal generator of $x_{n+2}J_n$. By cyclic symmetry of $C_n$
    we can assume that the element of  $x_{n+1}I_n$ is $x_1x_2x_{n+1}$. The remaining factor $x_3x_4\cdots x_nx_{n+2}$  is not in  $x_{n+2}J_n$ because $x_3x_4\cdots x_n$
    does not contain one of the variables $\{x_1,x_2\}$ (by definition of the cover ideal). 

Now consider $f^2=x_1^2x_2^2\cdots x_n^2x_{n+1}^2x_{n+2}^2$.
    If $n$ is even this can be written as $  (x_1x_3\cdots x_{n-1}x_{n+2})(x_2x_4\cdots x_{n}x_{n+2})
    (x_1x_2x_{n+1})(x_3x_4x_{n+1}) (x_5\cdots x_n) \in L_n^4,$  and  so $f\in \overline{L_n^2}\backslash L_n^2$. If $n$ is odd, then we have $   f^2=$
$$= (x_1x_2x_{n+1})(x_1x_3\cdots x_{n}x_{n+2})(x_2x_4\cdots x_{n-1}x_nx_{n+2})(x_3x_4x_{n+1})(x_5\cdots x_{n-1})\in$$ $\in L_n^4,$
     and so again $f\in \overline{L_n^2}\backslash L_n^2$.
This finishes our proof.
\end{proof}

\noindent  As an easy consequence we have the following result.

\begin{cor} If $n\geq 4$,  then $L_n$ is not normal.

\end{cor}

$\vspace{.1cm}$
\subsection{The case of $L_n^i, i\geq 3$}
$\vspace{.1cm}$

 In this subsection, we  prove the following result.

\begin{thm}\label{non-normal-1} If $n\geq 4,\ n$ even, or $n\geq 7,\ n$ odd, then $L_n^i, i\geq 3$ is not integrally closed.
\end{thm}

\begin{proof} 
We  produce an explicit element in  $\overline{L_n^i}\backslash L_n^i$ for $i\geq 3$ in these ranges of $n$.
  
 For $n$ even, $n\geq 4$, consider $h_i=x_1x_2^{i-1}x_3x_4^{i-1}\cdots x_{n-1}x_n^{i-1}x_{n+1}x_{n+2}^{i-1}$.
    Because of the factor $x_{n+1}x_{n+2}^{i-1}$,  $h_i$ is potentially in 
      $(x_{n+1}I_n)(x_{n+2}^{i-1}J_n^{i-1})\subset L_n^i$. 
    By cyclic symmetry we can assume that the
    generator in $x_{n+1}I_n$ is $x_1x_2x_{n+1}$. Dividing by this we get $x_2^{i-2}x_3x_4^{i-1}\cdots  x_{n-1}x_n^{i-1}x_{n+2}^{i-1}$. The latter cannot belong to $x_{n+2}^{i-1}J_n^{n+1}$, since we need $i-1$ factors of $x_1$ or $x_2$.
      Thus $h_i\notin L_n^i$. But we can write $$h_i^2=x_1^2x_2^{2i-2}x_3^2x_4^{2i-2}\cdots x_{n-1}^2x_n^{2i-2}x_{n+1}^2x_{n+2}^{2i-2}=(x_2x_4\cdots x_nx_{n+2})^{2i-4}f^2\in L_n^{2i},$$
      since we saw in the proof of Theorem \ref{non-normal} that $f^2\in L_n^4$. Thus $h_i\in \overline{L_n^i}\backslash L_n^i$, as desired.

      For $n\geq 7$ odd, let $g_i=x_2x_3x_4\cdots x_{n-2}x_{n-1}^{i-1}x_{n}^{i-1}x_{n+1}^{i-1}x_{n+2}$. If  $g_i\in L_n^i$ then
      $g_i\in (x_{n+1}I_n)^{i-1}(x_{n+2}J_n)$, because of the factor $x_{n+1}^{i-1}x_{n+2}$. If $p\in \mathcal G(x_{n+2}J_n)$, then $p$ must be divisible by at least one generator from each
      of the pairs $\{x_1,x_2\}$,$\{x_2,x_3\}$, $\ldots$, $\{x_{n-1}$,$x_n\}, \{x_n,x_1\}$. Since $x_1\notin$ supp($g_i$), $p$ must be divisible by $x_2$ and $x_n$.
      Furthermore $g_i/p$ cannot be divisible by any of the products $x_1x_2, x_2x_3,\ldots, x_{n-3}x_{n-2}$. The only part of $g_i/p$ left to give an element of $(x_{n+1}I_n)^{i-1}$
      is possibly $x_{n-2}x_{n-1}^{i-1}x_n^{i-2}x_{n+1}^{i-1}$. Here the exponents of $x_n$ and $x_{n+1}$ are certain but those of $x_{n-2}$ and $x_{n-1}$ are possibly one lower.
      Thus $g_i/p$ is always divisible by $(x_{n-1}x_n)^{i-2}$. Dividing by this we are left at most with $x_{n-2}x_{n-1}$. But $p$ has to be divisible by at least one
      of these variables, so $g_i/p$ cannot be in  $(x_{n+1}I_n)^{i-1}$, and $g_i$ cannot be in  $(x_{n+1}I_n)^{i-1}(x_{n+2}J_n)$ and hence is not in $L_n^i$.
      
      Note that $g_i$ is of degree $n-2+3(i-1)$ and that $x_{n+2}J_n$ has minimal generators of
      degree $2+(n-1)/2$ (e.g. $x_1x_3x_5x_7x_9$ for $C_7$). Thus potentially $g_i^2$ is the product of two minimal generators of $x_{n+2}J_n$ and $6(i-1)$ minimal generators of
      $x_{i+1}I_n$, with $2n-4+6(i-1)-(4+(n-1)+6(i-1))=n-7$ extra variables. Indeed, 
      \begin{align*}
      g_i^2=&(x_2x_3x_5x_7\cdots x_nx_{n+2})(x_2x_4x_6\cdots x_{n-1}x_nx_{n+2})(x_3x_4x_{n+1})(x_5x_6x_{n+1})\\
    &      (x_{n-1}x_nx_{n+1})^{2i-4}(x_7\cdots x_{n-1}).
      \end{align*}
      Thus, $g_i^2\in (x_{n+2}J_n)^2(x_{n+1}I_n)^{2i-2}\subset  L^{2i}$, and so $g_i\in \overline{L_n^i}$. The final conclusion is that $g_i\in \overline{L_n}^i\backslash L_n^i$.
     \end{proof}  

The last part of the proof above is not valid for $n=5$. Therefore we tackle this case in the next  subsection using the Rees cone.

\subsection{The Rees Cone, illustrated by application to $L_5$}
$\vspace{.1cm}$

 The Hilbert basis of the Rees cone can be computed with Normaliz \cite{N} (as in Example 2.16 of version 3.8.4).  
  A definition of the Rees cone and its properties can be found in \cite[Chapters 13,14]{VI}. 
  Let $I$ be a monomial ideal in a polynomial ring 
  ${\bQ}[x_1,\ldots,x_{n+2}]$
  whose minimal generators have exponent vectors $\{v_1,\ldots,v_m\}$. Let $\{e_i, 1\leq i\leq n+3\}$ be the unit vectors in ${\bQ}^{n+3}$. Then the Rees cone $C(I)$ is the rational cone in ${\bQ}^{n+3}$ spanned by the vectors $(v_i,1), 1\leq i\leq m$  and $e_i, 1\leq i\leq n+2$. This has the property that $(v,d)\in C(I)$ if and only if $\mathbf{x}^v\in\overline{I^d}$.
  The next example is motivated by Remark 1.3 (iii) of \cite{ANR}. 
  
  Using Normaliz directly we find that the Hilbert basis of the Rees cone of $I_5$
  is the rows of the matrix

\begin{verbatim}
       0 0 0 0 0 0 1 0
       0 0 0 0 0 1 0 0
       0 0 0 0 1 0 0 0
       0 0 0 1 0 0 0 0
       0 0 0 1 1 1 0 1
       0 0 1 0 0 0 0 0
       0 0 1 1 0 1 0 1
       0 1 0 0 0 0 0 0
       0 1 0 1 1 0 1 1
       0 1 1 0 0 1 0 1
       0 1 1 0 1 0 1 1
       1 0 0 0 0 0 0 0
       1 0 0 0 1 1 0 1
       1 0 1 0 1 0 1 1
       1 0 1 1 0 0 1 1
       1 1 0 0 0 1 0 1
       1 1 0 1 0 0 1 1
       1 1 1 1 1 1 1 2
\end{verbatim}

There are four groups of rows in this matrix. First, we have the $e_i$ for $1\leq i\leq 7$. Then there are the $(v_j, 1)$ for $1\leq j\leq 5$, where $v_j$ is the exponent vector of a minimal generator of $x_6I_5$, e.g. (0, 0, 0, 1, 1, 1, 0, 1) corresponds to $x_4x_5x_6$. There are also five exponent vectors of the minimal generators of $x_7J_5$, e.g. (0, 1, 0, 1, 1, 0, 1, 1) corresponds to $x_2x_4x_5x_7$. 
The final row corresponds to $x_1x_2x_3x_4x_5x_6x_7$, which according to the theory
of the Rees cone (and Theorem \ref{non-normal} above) is in  $\overline{L_5^2}\backslash L_5^2$.

\begin{thm}\label{normal} If $i \geq 3$,  then $L_5^i$ is integrally closed.
\end{thm}

\begin{proof} What we must show is that if $(v,i), i\geq 3$ is in the Rees cone of $I_5$, then $(v,i)$ is the sum of $i$ rows of the above matrix
  that end in 1 (possibly with some rows ending in 0). If there is no row $(1, 1, 1, 1, 1, 1, 1, 2)$, then we are done.
  At most one $(1, 1, 1, 1, 1,1,1,2)$ is required in this summation, by the proof of Theorem \ref{non-normal}.
  Thus, it suffices to show that $(1,1, 1, 1, 1, 1, 1, 2)$ plus any row ending in 1 can be rewritten without the $(1, 1, 1, 1, 1, 1, 1, 2)$. By cyclic symmetry
  in the variables $\{x_1,\ldots,x_5\}$ is suffices to consider one minimal generator of $x_6I_5$ and one of $x_7J_5$. Thus we   have 
  \begin{align*}
  (1, 1, 1, 1, 1, 1, 1, 2)+(0,0,0,1,1,1,0,1)=
  (1,1,1,2,2,2,1,3)=\\
 = (0,1,0,1,1,0,1,1)+(1,0,0,0,1,1,0,1) +  (0,0,1,1,0,1,0,1),
   \end{align*}
      and 
      \begin{align*}
      (1, 1, 1, 1, 1, 1, 1, 2)+(0, 1, 0, 1, 1, 0, 1, 1)= (1,2,1,2,2,1,2,3)=\\
 =(1,1,0,1,0,0,1,1)+(0,1,1,0,1,0,1,1)   + (0,0,0,1,1,1,0,1).
    \end{align*} 
\end{proof}

We have already seen  that  $f=x_1x_2x_3x_4x_5x_6x_7$ is not in $L_5^2$. In fact we have the following stronger result. The proof works for any odd $n\geq 5$, 
so we state it in that generality.

\begin{thm}\label{non-normal-3} Let $f=x_1x_2\cdots x_nx_{n+1}x_{n+2}$ where $n$ is an odd integer $\geq 5$. 
  For any $i\geq 1, x_{n+2}^if\notin L_n^2$, so that $\overline{L_n^2}/L_n^2$ is an infinite dimensional vector space.
\end{thm}

\begin{proof} We have $x_{n+2}^if=x_1x_2\cdots x_{n+1}x_{n+2}^{i+1}$. If this is in $L_n^2$ then from the $x_{n+1}x_{n+2}^{i+1}$ portion it must be in either
  $(x_{n+2}J_n)^2$ or $x_{n+1}I_nx_{n+2}J_n$. It cannot  be in the former because all minimal generators of $J_n$ are of degree $\geq (n+1)/2$, and $x_{n+2}^if$ is only of degree $n$ in $\{x_1,\ldots,x_n\}$.
  It cannot be in $x_{n+1}I_nx_{n+2}J_n$ by the same argument that $f\notin L_5^2$ (proof of Theorem \ref{non-normal}). In the proof of Theorem \ref{non-normal} it was seen that
  $f\in\overline{L_n^2}$ and of course $x_{n+2}^if$  is still in $\overline{L_n^2}$. Finally note that all ideals involved here are monomial ideals, and those monomials in $\overline{L_n^2}$
  but not in $L_n^2$ form a ${\bQ}$-basis of  $\overline{L_n^2}/L_n^2$. We have found an infinite number of these, so the conclusion follows.
\end{proof}  

We also have

\begin{thm}\label{non-normal-3a} $L_5^2:\overline{L_5^2}=(x_1,\ldots,x_6).$
\end{thm}

\begin{proof} We have $x_1f=x_1^2x_2x_3x_4x_5x_6x_7=(x_1x_2x_4x_7)(x_1x_5x_6)x_3\in L_5^2$. By cyclic symmetry in $\{x_1,\ldots,  x_5\}$ we also have
  $x_2f, x_3f, x_4f, x_5f \in L_5^2$. We also have $$x_6f=x_1x_2x_3x_4x_5x_6^2x_7=(x_1x_5x_6)(x_2x_3x_6)x_4x_7\in L_5^2.$$ 
  As in the proof of Theorem \ref{normal},  any monomial in $L_5^2$ can be written with at most one factor $f$.
  It now follows that $(x_1,\ldots,x_6)\in L_5^2:\overline{L_5^2}$. The only monomial ideals between  $(x_1,\ldots,x_6)$ and $(1)$ are of the form
  $(x_1,\ldots,x_6,x_7^i),$ $i\geq 1$. By Theorem \ref{non-normal-3}, $x_7^i\notin L_5^2:\overline{L_5^2}$. The claim now follows.
\end{proof}

It should be noted that, by  using similar techniques, we have been able to prove that $L_4^i:\overline{L_4^i}=(x_1,\ldots,x_6)$ for all  $i\geq 2$ and $\overline{L^i_4 }\setminus L^i_4$ is a finite dimensional vector  space for all $i\geq 1$. 

\noindent{\textbf{Acknowledgments.}}
The authors are deeply grateful to the  referee for careful reading of the manuscript and    valuable suggestions which improved  the quality of this paper.


\end{document}